# Easy representation of multivariate functions with low-dimensional terms via Gaussian process regression kernel design: applications to machine learning of potential energy surfaces and kinetic energy densities from sparse data

Eita Sasaki, Manabu Ihara and Sergei Manzhos *

School of Materials and Chemical Technology, Tokyo Institute of Technology, Ookayama 2-12-1, Meguro-ku, Tokyo 152-8552 Japan

*E-mail: E-mail: manzhos.s.aa@m.titech.ac.jp , Tel & Fax : 03-5734-3918

**We show that Gaussian process regression (GPR) allows representing multivariate functions with low-dimensional terms via kernel design. When using a kernel built with HDMR (High-dimensional model representation), one obtains a similar type of representation as the previously proposed HDMR-GPR scheme while being faster and simpler to use. We tested the approach on cases where highly accurate machine learning is required from sparse data by fitting potential energy surfaces and kinetic energy densities.**

Representations of a multivariate function $f(x), x \in R^D$ with low-dimensional terms are advantageous, as low-dimensional functions are easier to construct, especially from sparse data, and are advantageous in certain applications, most notably when the function needs to be integrated, such as in quantum dynamics calculations.[1,2] One approach for such a representation is High-dimensional model representation (HDMR) formalized by Rabitz and co-workers,[3–5] which is constructed as a sum of terms depending on subsets of original coordinates $(x_{i_1}, x_{i_2}, \ldots, x_{i_d}), d < D$:

$$f(x) \approx f_0 + \sum_{i=1}^{D} f_i(x_i) + \sum_{1 \leq i < j \leq D} f_{ij}(x_i, x_j) + \cdots + \sum_{\{i_1 i_2 \ldots i_d\} \in \{12 \ldots D\}} f_{i_1 i_2 \ldots i_d}(x_{i_1}, x_{i_2}, \ldots, x_{i_d})$$

(1)

We specifically consider RS (random sampling) HDMR[3,5] which allows obtaining all the terms from one and the same dataset with an arbitrary distribution of data in the full $D$-dimensional space. In the original RS-HDMR formulation, the component functions $f_{i_1 i_2 \ldots i_d}$ are obtained as $D$-$d$ dimensional integrals, which rapidly become a bottleneck as $D$ increases.[5–8] We previously introduced combinations of HDMR with neural networks (RS-HDMR-NN)[9–11] and, recently, Gaussian process regressions (RS-HDMR-GPR)[12,13] which allow dispensing with integrals and also allow combining terms of any dimensionality, e.g. one may use

$$f(x) \approx \sum_{\{i_1 i_2 \ldots i_d\} \in \{12 \ldots D\}} f^{NN/GPR}_{i_1 i_2 \ldots i_d}(x_{i_1}, x_{i_2}, \ldots, x_{i_d})$$

(2)

where the component functions are fitted with NN or GPR one at a time in cycles until convergence to the known values of the function $f^{(j)} = f(\pmb{x}^{(j)})$ at points $\pmb{x}^{(j)}, j = 1, …, M$:

$$f_{k_1 k_2 … k_d}(x_{k_1}, x_{k_2}, …, x_{k_d}) = f(\pmb{x}) - \sum_{\substack{\{i_1 i_2 … i_d\} \in \{12…D\} \\ \{i_1 i_2 … i_d\} \neq \{k_1 k_2 … k_d\}}} a(c) f_{i_1 i_2 … i_d}(x_{i_1}, x_{i_2}, …, x_{i_d})$$

(3)

where a fading parameter $a(c)$ (where $c$ indexes fitting cycles) may be introduced to palliate local minima.[12,13] In this way, existing NN or GPR engines can easily be used for component functions, although a custom code is needed to realize Eq. (3).[10,13] In eq. (2), lower order terms of eq. (1) are lumped into $d$-dimensional terms. This representation is particularly attractive with sparse data; we previously showed that with low data density, there is a maximum $d$ for which $f_{i_1 i_2 … i_d}$ can be reliably constructed.[9] Note that data density is always low in sufficiently high-dimensional spaces in any practical setting, as due to the so-called "curse of dimensionality" adding more data (even if it were possible, which is often not the case) has little effect on data density in terms of number of data points per dimension.[14] The representation of eq. (2) was previously used to fit functions to data with densities as low as about 2 data per degree of freedom or leess.[12,15] Eq. (3) does not enforce orthogonality of component functions but much gains in simplicity. When the terms are built with GPR (RS-HDMR-GPR[12,13]), eq. (2) allows determining relative importance of different combinations of variables, effectively extending the ARD (automated relevance determination) capability of plain GPR and allowing for pruning of HDMR terms;[13] this is important as the number of terms scales combinatorically with $D$ and $d$.

In eqs. (2-3), the relative amplitudes of terms are subsumed in the definition of $f_{i_1 i_2 … i_d}$; in what follows, it will be convenient to consider them explicitly:

$$f(\pmb{x}) \approx \sum_{\{i_1 i_2 … i_d\} \in \{12…D\}} A_{\widetilde{i_1 i_2 … i_d}} f_{\widetilde{i_1 i_2 … i_d}}(x_{i_1}, x_{i_2}, …, x_{i_d})$$

(4)

where $f_{\widetilde{i_1 i_2 … i_d}}$ are considered to be in some sense normalized (e.g. to have the maximum value or integral of 1). The amplitudes are fitted with eq. (3). The disadvantage of using eq. (3) is the need for repeated fits as well as the need for a separate code implementing the method. Another disadvantage is loss of ease of computing the variance of the estimate of $f(\pmb{x})$ (see eq. (6) below) which needs to be assembled from variance estimates of all component functions.

In GPR, the expectation values $f(\pmb{x})$ and variances $\Delta f(\pmb{x})$ of function values at any point in space $\pmb{x}$ are computed as[16]

$$f(\pmb{x}) = \pmb{K}^* \pmb{K}^{-1} \pmb{f} \tag{5}$$
$$\Delta f(\pmb{x}) = K^{**} - \pmb{K}^* \pmb{K}^{-1} \pmb{K}^{*T} \tag{6}$$

where $\pmb{f}$ is a vector of all (known) $f^{(j)}$ values, and the matrix $\pmb{K}$ and row vector $\pmb{K}^*$ are computed from pairwise covariances among the data:

$$K = \begin{pmatrix} k(x^{(1)}, x^{(1)}) + \delta & k(x^{(1)}, x^{(2)}) & \cdots & k(x^{(1)}, x^{(M)}) \\ k(x^{(2)}, x^{(1)}) & k(x^{(2)}, x^{(2)}) + \delta & \cdots & k(x^{(2)}, x^{(M)}) \\ \vdots & & \ddots & \vdots \\ k(x^{(M)}, x^{(1)}) & k(x^{(M)}, x^{(2)}) & \cdots & k(x^{(M)}, x^{(M)}) + \delta \end{pmatrix}$$

(7)

$$K^* = \begin{pmatrix} k(x, x^{(1)}) & k(x, x^{(2)}) & \cdots & k(x, x^{(M)}) \end{pmatrix},$$

(8)

and $K^{**} = k(x, x)$. The covariance function $k(x^{(1)}, x^{(2)}|\lambda)$ is the kernel of GPR that depends on hyperparameters $\lambda$ (which we omit in the formulas for notational simplicity). The optional $\delta$ on the diagonal has the meaning of the magnitude of Gaussian noise and is a regularization (hyper)parameter; it helps generalization.

Representation in the form of Eqs. (1-2) can also be obtained by GPR kernel design. Even though GPR is often considered to be a nonlinear machine learning method, at each particular value of $\lambda$, it is equivalent to a regularized linear regression. Eq. (5) has the form of a basis expansion,

$$f(x) = \sum_{n=1}^{M} b_n(x) c_n$$

(9)

with basis functions $b_n(x) = k(x, x^{(n)})$ and with coefficients $c$ obtained with least squares, $c = K^{-1} f$.[17] The covariance function is usually chosen as one of the Matern family of functions,[18]

$$k(x, x') = A \frac{2^{1-\nu}}{\Gamma(\nu)} \left( \sqrt{2\nu} \frac{|x - x'|}{l} \right)^\nu K_\nu \left( \sqrt{2\nu} \frac{|x - x'|}{l} \right)$$

(10)

where $\Gamma$ is the gamma function, and $K_\nu$ is the modified Bessel function of the second kind. At different values of $\nu$, this function becomes a squared exponential ($\nu \to \infty$), a simple exponential ($\nu = 1/2$) and various other widely used kernels (such as Matern3/2 and Matern5/2 for $\nu = 3/2$ and $5/2$, respectively). It is typically preset by the choice of the kernel, and the length scale $l$ and prefactor $A$ are hyperparameters (i.e. $\lambda = (l, A)$) that can be optimized. The particular value of GPR compared to a linear regression with a generic basis (which could in principle also be taken in an HDMR form) is the use of the covariance function which imparts Eqs. (5, 6) the meaning of the expectation value and the variance of a Gaussian distribution of $f(x)$ values.[16]

One can express $k$ in the form of eq. (1) or (2), e.g. in the case of eq. (2),

$$k(x, x') = \sum_{\{i_1 i_2 \ldots i_d\} \in \{12 \ldots D\}} A_{i_1 i_2 \ldots i_d} k_{i_1 i_2 \ldots i_d} (x_{i_1 i_2 \ldots i_d}, x'_{i_1 i_2 \ldots i_d})$$

(11)

where $x_{i_1 i_2 \ldots i_d} = (x_{i_1}, x_{i_2}, \ldots, x_{i_d})$ and $k_{i_1 i_2 \ldots i_d}$ can be chosen as one of Matern kernels in $d$ dimensions with amplitudes ($A$ in eq. (10)) indicated explicitly assuming $\max k_{i_1 i_2 \ldots i_d} = k_{i_1 i_2 \ldots i_d}(x_{i_1 i_2 \ldots i_d}, x_{i_1 i_2 \ldots i_d}) = 1$. Eq. (11) has previously been introduced for Additive Gaussian Processes.[19] Eqs. (9) and (11) together immediately give an HDMR-type representation of $f(x)$:

$$f(x) = \sum_{\{i_1 i_2 \ldots i_d\} \in \{12 \ldots D\}} A_{i_1 i_2 \ldots i_d} \sum_{n=1}^{M} k_{i_1 i_2 \ldots i_d}\left(x_{i_1 i_2 \ldots i_d}, x^{(n)}_{i_1 i_2 \ldots i_d}\right) c_n$$

$$\equiv \sum_{\{i_1 i_2 \ldots i_d\} \in \{12 \ldots D\}} B_{i_1 i_2 \ldots i_d} \tilde{f}_{i_1 i_2 \ldots i_d}(x_{i_1 i_2 \ldots i_d})$$

(12)

Note that elements of $c$, and thereby the amplitudes $B_{i_1 i_2 \ldots i_d}$, depend on products of many $A_{i_1 i_2 \ldots i_d}$ by virtue of the minors of the matrix $K$ forming its inverse (the dependence on $A_{i_1 i_2 \ldots i_d}$ of det($K$) need not be considered as it leads to a common scaling of all HDMR terms). That is, $B_{i_1 i_2 \ldots i_d} \neq A_{i_1 i_2 \ldots i_d}$, and explicit dependence of $B_{i_1 i_2 \ldots i_d}$ on $A_{i_1 i_2 \ldots i_d}$ is impractically complex. However, *the choice of $A_{i_1 i_2 \ldots i_d}$ is immaterial*. One can even choose the amplitudes of HDMR terms *of the kernel* randomly. Regardless of the choice of $A_{i_1 i_2 \ldots i_d}$, $B_{i_1 i_2 \ldots i_d}$ are the least squares solutions and are in this sense optimal (one can think of any changes introduced in relative magnitudes of different $A_{i_1 i_2 \ldots i_d}$ being compensated in Eq. (12) via $c$). This is in contrast to Eq. (3) where the amplitudes of the component functions depend on the quality of the optimization and local minima. Eq. (11) as written (which is the form we use in the tests below) uses all combinations of $d$ variables. The approach is obviously not limited to this particular form; a generic HDMR expansion of the form of eq. (1) can also be used for the kernel and will result in a corresponding HDMR expansion of $f(x)$. Only selected combinations of $d' \leq d$ variables can also be used to decrease the number of terms.[13] Individual component functions are computable as

$$f_{k_1 k_2 \ldots k_d}(x_{k_1}, x_{k_2}, \ldots, x_{k_d}) = K^*_{i_1 i_2 \ldots i_d} c$$

(13)

where $K^*_{i_1 i_2 \ldots i_d}$ is a row vector with elements $A_{i_1 i_2 \ldots i_d} k_{i_1 i_2 \ldots i_d}\left(x_{i_1 i_2 \ldots i_d}, x^{(n)}_{i_1 i_2 \ldots i_d}\right)$. In particular, the values of the component functions at the training set are $f_{i_1 i_2 \ldots i_d} = K_{i_1 i_2 \ldots i_d} c$ and can be used to evaluate the relative importance of different component functions by computing the variance $var(K_{i_1 i_2 \ldots i_d} c)$, where the (m,n) elements of the matrix $K_{i_1 i_2 \ldots i_d}$ are $A_{i_1 i_2 \ldots i_d} k_{i_1 i_2 \ldots i_d}\left(x^{(m)}_{i_1 i_2 \ldots i_d}, x^{(n)}_{i_1 i_2 \ldots i_d}\right)$.

Rather than using a dedicated code as in the case of Eq. (3),[13] an HDMR-type kernel of Eq. (11) can be used with any existing GPR engine to obtain a HDMR representation of $f(x)$ directly. One simply needs to define a custom kernel function, which is easily doable in common machine learning libraries such as Matlab's Statistics and Machine Learning Toolbox used by us. The use of a single GPR approximation instead of Eq. (3) also makes easy the calculation of the variance of $f(x)$ with Eq. (6), which is also then returned by the GPR engine. We caution that Eq. (6) should not be automatically used as an error bar; it allows computing the confidence interval on the expectation value of the function computed by GPR and might not be indicative of the quality of the fit (see notes and examples to this effect in Refs. [12,13]).

Table 1. Test set rmse errors obtained with RS-HDMR-GPR (Eq. (3)) of different orders $d$ in Refs. [12,13] for the potential energy surfaces of $H_2CO$ and $UF_6$ and for the kinetic energy densities of Al, Mg, and Si as well as those obtained with an HDMR-type kernel in this work for different numbers of training datapoints $M$. The $M$ in the case of RS-HDMR-GPR results are the largest among those used in Refs. [12,13]. Test set sizes are 100,000 for $H_2CO$ PES, 40,000 for $UF_6$ PES, and 400,000 for the KEDs. $N$ is the number of HDMR component functions at each $d$. The numbers are in cm$^{-1}$ for the PESs and in *a.u.* for KED. The ranges of rmse from 10 runs are given (where available) reflecting the random nature of the training data selection from the overall dataset.

| | | $H_2CO$ PES ($D = 6$) | | | | $UF_6$ PES ($D = 15$) | | | | KED ($D = 7$), in units of ×10$^{-4}$ | | |
|---|---|---|---|---|---|---|---|---|---|---|---|---|
| $d$ | $N$ | Ref. [12] | Eq. (11) | Eq. (11) | $N$ | Ref. [12] | Eq. (11) | Eq. (11) | $N$ | Ref. [13] [a] | Eq. (11) | Eq. (11) |
| $M$ | | 3600 | 3600 | 10000 | | 5000 | 5000 | 10000 | | 5000 | 5000 | 10000 |
| $d = D$ | 1 | 1.58-1.79 | 0.99-1.23 | 0.37-0.56 | 1 | 42.2 | 36.0-39.2 | 25.5-26.2 | 1 | 2.53 | 2.34-6.47 | 2.21-2.44 |
| 1 | 6 | 1277-1280 | 1274-1280 | 1271-1276 | 15 | 234.6 | 234-236 | 234-235 | 7 | 9.10 | 4.84-4.86 | 4.82-4.85 |
| 2 | 15 | 295-300 | 285-289 | 276-279 | 105 | 168.1 | 165-168 | 161-163 | 21 | 4.46 | 2.70-4.32 | 2.61-2.97 |
| 3 | 20 | 10.70-11.25 | 8.25-10.6 | 8.81-9.98 | 455 | 65.6 | 60.5-61.7 | 55.5-56.3 | 35 | 2.72 | 2.57-3.33 | 2.47-2.68 |
| 4 | 15 | 1.36-1.83 | 0.74-1.05 | 0.32-0.45 | 1365 | | 58.6-59.5 | 51.8-52.5 | 35 | 2.72 | 2.48-2.69 | 2.38-2.87 |
| 5 | 6 | 1.23-1.78 | 0.94-1.16 | 0.34-0.40 | | | | | 21 | 1.73 | 2.45-3.24 | 2.32-2.65 |
| 6 | | | | | | | | | 7 | 1.73 | 2.38-3.45 | 2.26-2.41 |
| … | | | | | | | | | | | | |
| 11 | | | | | 1365 | | 36.1-37.3 | 25.5-25.9 | | | | |
| 12 | | | | | 455 | | 36.1-37.0 | 25.3-25.8 | | | | |
| 13 | | | | | 105 | | 35.7-37.1 | 25.5-25.9 | | | | |
| 14 | | | | | 15 | | 35.8-38.2 | 25.5-26.1 | | | | |

We compared the performance of an HDMR-type kernel to that of the RS-HDMR-GPR[12,13] approach. We fitted the potential energy surfaces (PES) of $H_2CO$ and $UF_6$ molecules and the kinetic energy densities (KED) of aluminium, magnesium, and silicon crystals at equilibrium and strained geometries. For the description of the datasets and information about the applications of these functions, see Refs. [20] for $H_2CO$, Ref. [21] for $UF_6$, and Ref. [22] for the KED data. We chose these applications, as in them, high fitting accuracy is required (with errors much smaller than 1% of the data range for the model to be useful at all and with desirable accuracy of better than 0.01%[23–25]). We have 120,000 data points in six dimensions (representing molecular bonds and angles) for $H_2CO$ with values of potential energy ranging 0-17,000 cm$^{-1}$, 54,991 data points in 15 dimensions (representing 15 modes of vibration) for $UF_6$ with values ranging 0-6,629 cm$^{-1}$, and 585,890 data in seven dimensions (representing six terms of the 4$^{th}$ order gradient expansion of kinetic energy density[26] and the effective potential[27]) for KED with values ranging 0.0073-0.0398 *a.u.* The distributions of the data can also be found in the original references. What matters for the purpose of the present work are not details of these applications but a comparison of GPR with HDMR-type kernel of Eq. (11) to RS-HDMR-GPR and plain GPR (obtained when $d = D$ in eq. (11)) using a standard kernel on the same data.

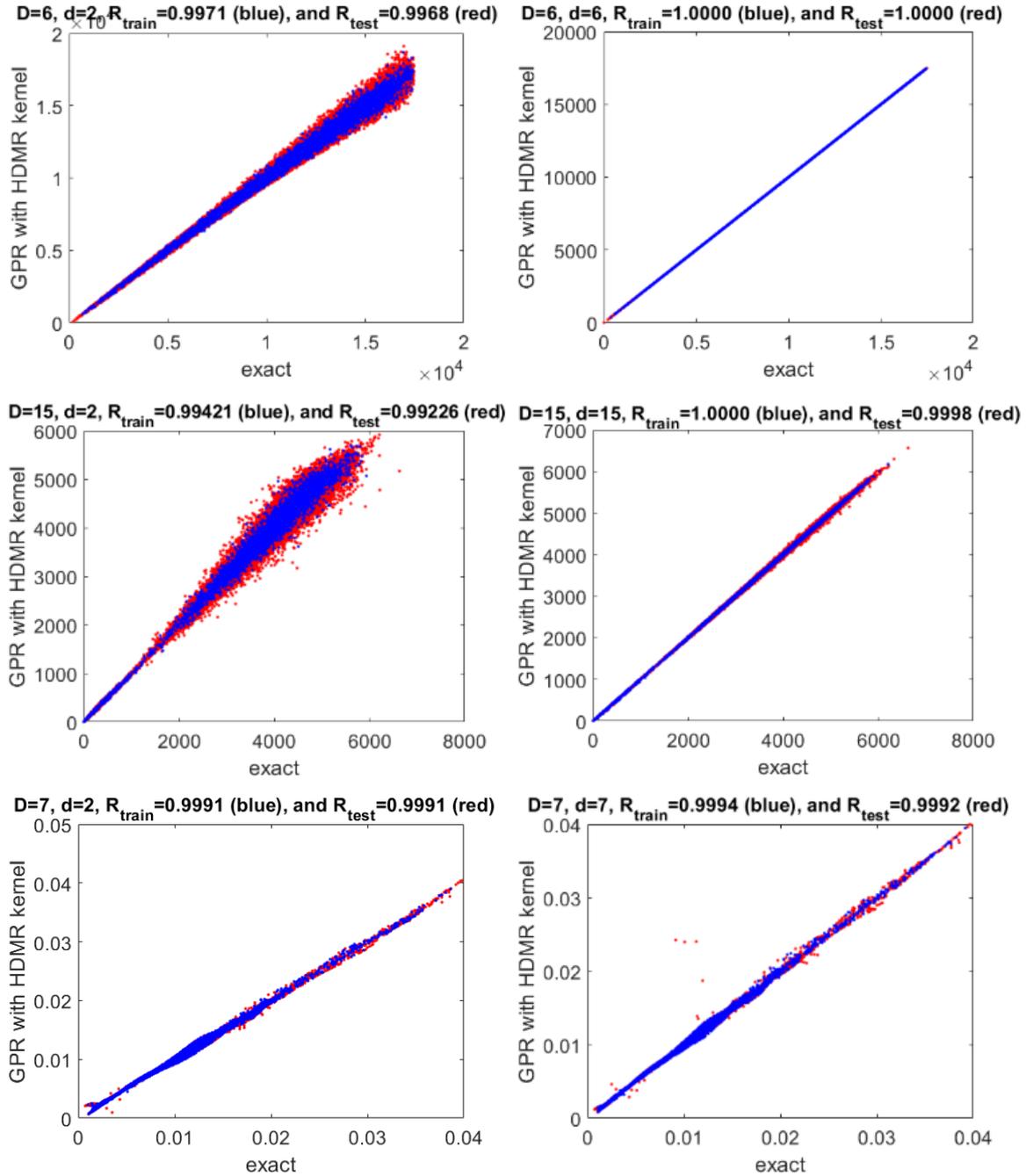

Figure 1. Representative examples of correlation plots and correlation coefficients between exact values and HDMR kernel based GPR predictions for selected combinations of *d* and *D* (for the case of 10,000 training points).

For $k_{i_1 i_2 \ldots i_d}(x_{i_1 i_2 \ldots i_d}, x'_{i_1 i_2 \ldots i_d})$ of Eq. (11), we use the squared exponential kernel (i.e. $k_{i_1 i_2 \ldots i_d}(x_{i_1 i_2 \ldots i_d}, x'_{i_1 i_2 \ldots i_d}) = exp\left(|x_{i_1 i_2 \ldots i_d} - x'_{i_1 i_2 \ldots i_d}|^2/2l^2\right)$). The PES data were normalized before fitting; we therefore use isotropic kernels. The KED data were scaled to [0,1] for the same reason;

as their distributions are extremely uneven (see Ref. [22]), scaling is preferred over normalization in this case. The length parameter $l$ is 5.47 ($d$ = 4-6) – 8.17 ($d$ = 1-3) for $H_2CO$, 33.1 for $UF_6$, and 1.22 for KED. The corresponding $\delta$ values are $1\times10^{-6}$ ($d$ = 4-6) - $1\times10^{-5}$ ($d$ = 1-3) for $H_2CO$, $1\times10^{-5}$ for $UF_6$, and $5\times10^{-4}$ for KED. We set all $A_{i_1 i_2 \ldots i_d}$ to the same value ($1/N$ where $N = C_d^D$ is the number of HDMR terms and $C_d^D$ is the binomial coefficient). We confirmed that they can be changed randomly without affecting the quality of the fit, just as the theory suggests. There is no effect of the setting of the relative amplitudes of component functions in the kernel; however, there is an effect via $\delta$ as the effect of $\delta$ depends on the overall magnitude of the kernel. Setting all $A_{i_1 i_2 \ldots i_d} = const$ is sufficient; *there is no need to optimize the magnitudes of the kernel's component functions*.

The tests were run in Matlab 2021a using *fitrgp* function with a custom kernel function implementing Eq. (11). In Table 1, we compare test set root mean square errors (rmse) obtained with the kernel of Eq. (11) to those reported in Refs. [12,13] with RS-HDMR-GPR (using similar kernels). We use test set sizes (also given in Table 1) which are much larger than the training sets to grasp well the global quality of the approximation. Note that there is variability of rmse values from run to run due to random selection of training data from the overall data set, which is within about ±10% and does not affect the reported trends (we do not fix the random seed precisely to monitor the effect of this source of uncertainly and give rmse ranges from 10 runs). For all practical purposes, the fit quality is similar to that achieved with RS-HDMR-GPR.[12,13] The obtainable error for given hyperparameters with the HDMR kernel is lower than with RS-HDMR-GPR by construction, as the coefficients are optimal in the least squares sense. In the case of the KED data, the errors obtained here for $d < D$ are slightly higher than those reported in Ref. [13], as in that work, length parameters were optimized for each component function (see below for tests with optimized $l$).

Table 2. Relative variances of the seven component functions of the HDMR obtained with the kernel of Eq. (11) for $d = 1$ when fitting 10,000 KED data with fixed and optimized length parameters. The values are in *a.u*.

| $l = 1.22$ | Optimized $l$ | |
|---|---|---|
| $var(K_i c) \times 10^4$ | $var(K_i c) \times 10^4$ | $l_{opt}$ |
| 0.00960 | 0.00897 | 0.297 |
| 0.00088 | 0.00082 | 0.099 |
| 1.03933 | 1.06364 | 0.135 |
| 0.00024 | 0.00001 | 6.058 |
| 0.00004 | 0.00000 | 21.93 |
| 0.00042 | 0.00046 | 0.693 |
| 1.73715 | 1.75074 | 0.131 |
| Train/test rmse $\times 10^4$ | | |
| 4.82 / 4.83 | 4.23 / 4.32 | |

When using an HDMR-type kernel, it is easier to use larger training sets, as a single GPR instance is fitted once. In Refs. [12,13], a maximum of 3600, 5000, and 5000 points were fitted for $H_2CO$, $UF_6$, and KED, respectively. Larger sets were not used in Refs. [12,13], in particular, due to the relatively high scaling of cost of GPR with the number of training data, compounded by the cost of applying Eq. (3) and wielding $N = C_d^D$ GPR instances. With an HDMR-type kernel, the cost is still higher than that of a conventional Matern-type kernel due to a higher cost of computing the kernel function which has a larger number of terms (the number of terms is also given in Table 1 for each $d$) but is easier manageable. We also provide results with much larger training sets and larger $d$ in the case of $UF_6$, where HDMR has more than a thousand terms. As expected, even lower global errors are obtainable for higher $d$ with more training data, while lower-dimensional terms are well-defined with few data[9] (i.e. the test rmse is not improved by adding more data and is limited by the dimensionality of HDMR terms). These results with the HDMR kernel highlight the advantages of the HDMR representation, namely, that with finite training data, one can obtain a similar or better global rmse compared to a conventional (full-dimensional) GPR with $d < D$. We show representative examples of correlation plots between exact values and HDMR kernel based GPR predictions for select $d < D$ in Figure 1, to visually highlight the high accuracy of the regressions performed here.

We mentioned above that the HDMR-type kernel allows estimating variances of component functions and therefore, similar to RS-HDMR-GPR, relative importance of different variables or combinations thereof. Taking the KED data as an example, we list the variances of the 7 component functions at $d = 1$ obtained at fixed and optimized (to maximum likelihood) length parameters in Table 2. Similar to what was found with RS-HDMR-GPR in Ref. [13], variables $x_1, x_2, x_3, x_7$ are seen as most important. Specifically, the dwindling of the variance of $f_4(x_4)$ and $f_5(x_5)$ corresponds to their optimized length parameters becoming large, indicating their low relevance, in a way similar to ARD.

## Conclusions

We showed that a kernel type based on High-dimensional model representation for Gaussian process regression allows easily constructing a representation of a multivariate function as a sum of lower-dimensional terms. A similar kernel representation was previously introduced for Additive Gaussian Processes;[19] here, we show that it allows obtaining accurate models in applications notorious for high accuracy requirements – potential energy surfaces and kinetic energy densities and allows building HDMR representations of multivariate functions which are similar to the recently proposed RS-HDMR-GPR scheme[12,13] while being much easier to use. One only needs to define a custom kernel for use with existing GPR libraries; no dedicated software is needed. There is no need to optimize the magnitudes of the HDMR terms in the kernel, and the magnitudes of the HDMR term of the final function representation are optimal in the least squares sense. The relative importance of different HDMR component functions and corresponding variables can also be easily evaluated.

## Data availability

Matlab codes and datasets used in this work are available from the authors upon request.

## Acknowledgements

We thank Prof. Tucker Carrington and Dr. Owen Ren for discussions. We thank Dr. Laura Laverdure and Prof. Nicholas Mosey who computed the data for Ref. [21].